\documentclass[reqno]{amsart}
\usepackage{amsmath,amssymb,amsthm,ascmac}
\usepackage[mathscr,mathcal]{eucal}
\usepackage{latexsym}
\usepackage[dvipdfmx]{graphicx}
\usepackage{enumitem}

\begin{document}

\allowdisplaybreaks

\theoremstyle{definition}
\newtheorem{defi}{\textbf{Definition}}[section]
\newtheorem{thm}[defi]{\textbf{Theorem}}
\newtheorem{lem}[defi]{\textbf{Lemma}}
\newtheorem{prop}[defi]{\textbf{Proposition}}
\newtheorem{cor}[defi]{\textbf{Corollary}}
\newtheorem{ex}[defi]{\textbf{Example}}
\newtheorem{rem}[defi]{\textbf{Remark}}
\newtheorem*{corr}{\textbf{Corollary}}

\theoremstyle{plain}
\newtheorem{maintheorem}{Theorem}
\newtheorem{theorem}{Theorem }[section]
\newtheorem{proposition}[theorem]{Proposition}
\newtheorem{mainproposition}{Proposition}
\newtheorem{lemma}[theorem]{Lemma}
\newtheorem{corollary}[theorem]{Corollary}
\newtheorem{maincorollary}{Corollary}
\newtheorem{claim}{Claim}
\renewcommand{\themaintheorem}{\Alph{maintheorem}}
\theoremstyle{definition} \theoremstyle{remark}
\newtheorem{remark}[theorem]{Remark}
\newtheorem{example}[theorem]{Example}
\newtheorem{definition}[theorem]{Definition}
\newtheorem{problem}{Problem}
\newtheorem{question}{Question}
\newtheorem{exercise}{Exercise}

\newtheorem*{subject}{�ړI}
\newtheorem*{mondai}{Problem}
\newtheorem{lastpf}{�ؖ�}
\newtheorem{lastpf1}{proof of theorem2.4}
\renewcommand{\thelastpf}{}

\newtheorem{last}{Theorem}

\renewcommand{\thelast}{}
\renewcommand{\proofname}{\textup{Proof.}}

\renewcommand{\theequation}{\arabic{section}.\arabic{equation}}
%%%%%  ���ԍ��� section ���Ƀ��Z�b�g����  %%%%%
\makeatletter
\@addtoreset{equation}{section}

\title[One-way specification and large deviations]
{On the one-way specification property and large deviations for non-transitive systems}
\author[K. Yamamoto]{Kenichiro Yamamoto}
\address{Department of General Education \\
Nagaoka University of Technology \\
Niigata 940-2188, JAPAN
}
\email{k\_yamamoto@vos.nagaokaut.ac.jp}


\date{}

\maketitle
\large

\begin{abstract}
We introduce a weaker form of the specification property, called ``one-way specification property", and give several examples of non-transitive systems
satisfying this property.
As an application, we show that the $(-\beta)$-transformation satisfies a level-2 large deviation principle with the Lebesgue measure
and the rate function is the free energy under the condition that
$\beta>1$ is a Yrrap number, that is, the orbit of $1$ under the $(-\beta)$-transformation
is eventually periodic.
\end{abstract}

\section{Introduction}

The specification property was introduced by Bowen in \cite{B} to study the ergodic property of Axiom A diffeomorphisms. After that several authors introduced
various weaker forms of this property such as weak specification property (\cite{D,M}), almost specification property (\cite{PfS2,Th}), specification property
for collection of good words (\cite{CT}) and
non-uniform specification property (\cite{V}) (see also \cite{KLO} for details).
These weaker forms of the property hold mainly for topologically transitive systems, for example, ergodic group automorphisms, $\beta$-shifts and $S$-gap shifts.
In this paper, we introduce a new weaker form of the specification property called ``one-way ((W)-) specification property", which holds for a large class of
non-transitive systems and investigate large deviations for systems with this property.

Let $(X,d)$ be a compact metric space and $f$ be  a Borel measurable map from $X$ to itself. We denote by $\mathcal{M}(X)$ the set of all Borel probability measures on $X$ with the weak topology. We say that $(X,f)$ satisfies a level-2 large deviation principle with a reference measure $m\in\mathcal{M}(X)$
if there is an upper semicontinuous function $q\colon\mathcal{M}(X)\to [-\infty,0]$, called a rate function, satisfying
\begin{displaymath}
\liminf_{n\rightarrow\infty}\frac{1}{n}\log m\left(\left \{x\in X:\frac{1}{n}\sum_{j=0}^{n-1}\delta_{f^j(x)}\in G\right \}\right)\ge\sup_{\mu\in G}q(\mu)
\end{displaymath}
holds for any open set $G\subset\mathcal{M}(X)$ and
\begin{displaymath}
\limsup_{n\rightarrow\infty}\frac{1}{n}\log m\left(\left \{x\in X:\frac{1}{n}\sum_{j=0}^{n-1}\delta_{f^j(x)}\in F\right \}\right)\le\sup_{\mu\in F}q(\mu)
\end{displaymath}
holds for any closed set $F\subset\mathcal{M}(X)$, where $\delta_y$ stands for the $\delta$-measure at the point $y\in X$.
For an $f$-invariant measure $\mu\in\mathcal{M}(X)$, let $h(\mu)$ be the metric entropy of $\mu$. For a continuous function
$\varphi\colon X\to\mathbb{R}$, we define a free energy function (with respect to $\varphi$),
$\mathcal{F}_{\varphi}\colon\mathcal{M}(X)\to [-\infty,0]$ by
\begin{equation}
\label{free}
\mathcal{F}_{\varphi}(\mu)=
\left\{
\begin{array}{ll}
h(\mu)-\int\varphi d\mu & (\mu\text{ is invariant}); \\
-\infty & (\text{otherwise}).
\end{array}
\right.
\end{equation}
Now we state our main theorem (see \S2.1 for precise definitions).
\begin{maintheorem}
\label{main1}
Let $(X,d)$ be a compact metric space, $f\colon X\to X$ be a continuous map, $m\in\mathcal{M}(X)$, $\varphi$ be a continuous function on $X$ and suppose that the four conditions \ref{A.I}-\ref{A.IIII} hold.
\noindent
\begin{enumerate}[label=\textbf{[A.\arabic{*}]}]
\item\label{A.I}
$(X,f)$ satisfies the one-way (W)-specification property.

\item\label{A.II}
$\displaystyle \lim_{\epsilon\rightarrow 0}\liminf_{n\rightarrow\infty}\inf_{x\in X}\left(\frac{1}{n}\log m(B_n(x,\epsilon))+\frac{1}{n}\sum_{j=0}^{n-1}\varphi(f^j(x))\right)\ge 0$.

\item\label{A.III}
The entropy map $\mu\mapsto h(\mu)$ is upper semicontinuous on the set of all $f$-invariant Borel probability measures on $X$.

\item\label{A.IIII}
$\displaystyle \lim_{\epsilon\rightarrow 0}\limsup_{n\rightarrow\infty}\sup_{x\in X}\left(\frac{1}{n}\log m(B_n(x,\epsilon))+\frac{1}{n}\sum_{j=0}^{n-1}\varphi(f^j(x))\right)\le 0$.
\end{enumerate}
Here $B_n(x,\epsilon):=\{y\in X:d(f^j(x),f^j(y))\le\epsilon,\ 0\le j\le n-1\}$. Then $(X,f)$ satisfies a level-2 large deviation principle with $m$ and
the rate function is the free energy $\mathcal{F}_{\varphi}\colon\mathcal{M}(X)\to [-\infty,0]$ given by (\ref{free}).
\end{maintheorem}

As a consequence of Theorem \ref{main1}, we have the following results for symbolic systems and piecewise expanding interval systems
(for precise definitions, see \S2.2 and \S2.3).

\begin{maincorollary}
\label{maincor}
Let $(X,\sigma)$ be a subshift on a finite alphabet, $m\in\mathcal{M}(X)$ and $\varphi\colon X\to \mathbb{R}$ be a continuous function.
Let $\mathcal{L}$ be the language of $X$.
Suppose that the three conditions
\ref{1.1}-\ref{1.3} hold.
\begin{enumerate}[label=\textbf{[1.\arabic{*}]}]
\item\label{1.1}
$(X,\sigma)$ satisfies the one-way (W)-specification property.

\item\label{1.2}
$\displaystyle\liminf_{n\rightarrow\infty}\inf_{w\in \mathcal{L},|w|=n}\left(\frac{1}{n}\log m([w])+\inf_{\omega\in [w]}\frac{1}{n}\sum_{j=0}^{n-1}\varphi(\sigma^j(\omega))\right)\ge 0$.

\item\label{1.3}
$\displaystyle\limsup_{n\rightarrow\infty}\sup_{w\in \mathcal{L},|w|=n}\left(\frac{1}{n}\log m([w])+\sup_{\omega\in [w]}\frac{1}{n}\sum_{j=0}^{n-1}\varphi(\sigma^j(\omega
))\right)\le 0$.
\end{enumerate}
Here $[w]$ denotes the cylinder set for $w\in\mathcal{L}$.
Then $(X,\sigma)$ satisfies a level-2 large deviation principle with $m$ and the rate function
is the free energy $\mathcal{F}_{\varphi}\colon\mathcal{M}(X)\to [-\infty,0]$ given by (\ref{free}).
\end{maincorollary}

\begin{maincorollary}
\label{p-e}
Let $([0,1],f)$ be a piecewise expanding interval system, $(\Sigma_f,\sigma)$ be a symbolic space of $([0,1],f)$,
$L$ be the Lebesgue measure on $[0,1]$ and $\varphi\colon [0,1]\to\mathbb{R}$ be a continuous function.
Suppose that the two conditions \ref{2.1} and \ref{2.2} hold.
\begin{enumerate}[label=\textbf{[2.\arabic{*}]}]
\item\label{2.1}
There exists a measurable map $\Phi\colon [0,1]\to \Sigma_f$ such that
$\Psi\circ\Phi\colon [0,1]\to [0,1]$ is an identity map and
$\Phi\circ f=\sigma\circ \Phi$.
Moreover, there exists a Borel set $Y\subset \Sigma_f$ such that
$[0,1]\setminus \Psi(Y)$ is countable and $\Phi\circ\Psi\colon Y\to \Phi(\Psi(Y))$ is an identity map.

\item\label{2.2}
$(\Sigma_f,\sigma)$ satisfies the conditions \ref{1.1}-\ref{1.3} of Corollary \ref{maincor} for $L_f:=L\circ\Phi^{-1}$
and $\varphi\circ\Psi$.
\end{enumerate}
Here $\Psi\colon \Sigma_f\to [0,1]$ is the natural factor map.
Then $(I,f)$ satisfies a level-2 large deviation principle with $L$
and the rate function is the free energy $\mathcal{F}_{\varphi}\colon\mathcal{M}([0,1])\to [-\infty,0]$ given by (\ref{free}).
\end{maincorollary}

In this paper, we give several examples of systems with the one-way ((W)-) specification property.
The main example is a $(-\beta)$-transformation. The $(-\beta)$-transformation was recently introduced by Ito and Sadahiro in \cite{IS}
and several authors considered problems of $(-\beta)$-transformations, which are also considered in $\beta$-transformations (\cite{DR,FL,LS}).
In the case of $\beta$-transformations, Pfister and Sullivan \cite{PfS} proved a level-2 large deviation principle.
As an application of Corollary \ref{p-e}, we prove the analogous result for $(-\beta)$-transformations
under the condition that $\beta>1$ is a Yrrap number.

\begin{maintheorem}
\label{minus beta}
Let $T_{-\beta}\colon [0,1]\to [0,1]$ be a $(-\beta)$-transformation, $L$ be the Lebesgue measure on $[0,1]$
and $\varphi\colon [0,1]\to \mathbb{R}$ be a constant function $\varphi\equiv \log \beta$.
Suppose that $\beta>1$ is a Yrrap number, that is, the orbit of $1$ under $T_{-\beta}$ is eventually periodic.
Then $([0,1],T_{-\beta})$ satisfies a level-2 large deviation principle
with $L$ and the rate function is the free energy $\mathcal{F}_{\varphi}\colon\mathcal{M}([0,1])\to [-\infty,0]$ given by (\ref{free}).
\end{maintheorem}
\begin{rem}
We note that our definition of the $(-\beta)$-transformation is slightly different from the original one \cite{IS} (see Remark \ref{rem-minus-beta}).
\end{rem}
Large deviations theory for dynamical systems has been studied by many authors and several authors proved a full large deviation principle such as
\cite{Ch,ChT1,CTY,EKW,PfS,T1,Ya,Y} (for other results, see \cite[\S1]{CTY}). The common property which holds for systems appeared in these papers are ``entropy-density", that is,
any invariant measure can be approximated by ergodic measures with similar entropies (see \cite{PfS} for the precise definition).
In the case of $(-\beta)$-transformations, if $\beta$ is a minimal Pisot number, then it is also a Yrrap number and the set of all $T_{-\beta}$-ergodic measures are not
entropy-dense (even not dense)
in the set of all invariant measures. Thus by Theorem \ref{minus beta}, $(-\beta)$-transformation for the minimal Pisot number $\beta$ satisfies a
level-2 large deviation principle with the Lebesgue measure and the rate function is the free energy though the set of ergodic measures are
not entropy-dense.
As far as we know, this is the first non-trivial example satisfying such properties.
On the other hand, as mentioned above, Pfister and Sullivan \cite{PfS} proved that the $\beta$-transformation ($\beta>1$) satisfies a level-2 large deviation principle
with a unique measure of maximal entropy and  the rate function is the free energy.
In this case, it is well-known that the unique measure of maximal entropy is
absolutely continuous with respect to Lebesgue measure, with a density bounded above by $\beta/(\beta-1)$ and below by $(\beta-1)/\beta$
(see \cite{P,R}).
Thus, the $\beta$-transformation also satisfies a large deviation principle with the Lebesgue measure and the rate function
is the free energy as in the case of the measure of maximal entropy. However, it occurs a different phenomenon in the case of $(-\beta)$-transformations.
Indeed, as a consequence of Theorem \ref{minus beta}, if $\beta$ is a minimal Pisot number, then
the $(-\beta)$-transformation satisfies a level-2 large deviation principle
with both the Lebesgue measure and the measure of maximal entropy, but these two rate functions do not coincide (see Example \ref{minimal pisot}).

This paper is organized as follows.
In \S2, we establish our definitions. In \S3, we give several examples of systems satisfying the one-way ((W)-) specification property.
In \S4, we give a proof of Theorem \ref{main1} and also give a proof of Corollaries \ref{maincor} and \ref{p-e} in \S5.
In \S6, we apply our results to systems appeared in \S3.1 and prove Theorem \ref{minus beta}.

\section{Preliminaries}
\subsection{Definitions, basic facts and lemmas}
Let $(X,d)$ be a compact metric space and $f\colon X\to X$ be a continuous map.
We denote by $C(X)$ the Banach space of continuous real-valued functions of $X$
with the sup norm $\|\cdot\|_{\infty}$ and by $\mathcal{M}(X)$ the set of all Borel probability measures on $X$ with
the weak topology. Since $C(X)$ is separable, there exists a countable set $\{\varphi_1,\varphi_2,\cdots\}$ which is dense in $C(X)$.
For $\mu,\nu\in\mathcal{M}(X)$, we define
$$D(\mu,\nu):=\sum_{n=1}^{\infty}\frac{|\int\varphi_n d\mu-\int\varphi_n d\nu|}{2^{n+1}\|\varphi_n\|_{\infty}}.$$
Then $D$ is a compatible metric for $\mathcal{M}(X)$.
By an easy calculation, we have the following.
\begin{lem}
\label{weak metric}
(1) $D(\mu,\nu)\le 1$ holds for any $\mu,\nu\in\mathcal{M}(X)$.

\noindent
(2) Let $\sum_{i=1}^pa_i\mu_i$ and $\sum_{i=1}^pa_i\nu_i$ be finite convex conbinations of Borel probability measures on $X$. Then we have
$$D\left(\sum_{i=1}^pa_i\mu_i,\sum_{i=1}^pa_i\nu_i\right)\le \sum_{i=1}^pa_iD(\mu_i,\nu_i).$$

\noindent
(3) Let $\sum_{i=1}^pa_i\mu_i$ and $\sum_{i=1}^pb_i\nu_i$ be finite convex conbinations of Borel probability measures on $X$.
Suppose that $\sum_{i=1}^p|a_i-b_i|\le \zeta$ and $D(\mu_i,\nu_i)\le \zeta$ for $1\le i\le p$.
Then we have
$$D\left(\sum_{i=1}^p a_i\mu_i,\sum_{i=1}^p b_i\nu_i\right)\le 2\zeta.$$
\end{lem}

Let $\mathcal{M}_f(X)\subset \mathcal{M}(X)$ be the set of all $f$-invariant
Borel probability measures, and let $\mathcal{M}_f^e(X)\subset\mathcal{M}_f(X)$
be the set of all ergodic measures.
For $n\ge 1$, we define $\mathcal{E}_n\colon X\to \mathcal{M}(X)$ by $\mathcal{E}_n(x):=\frac{1}{n}\sum_{j=0}^{n-1}\delta_{f^j(x)}$.

Before giving a definition of the one-way specification property, we recall the definition of the classical specification property.
We say that $(X,f)$ satisfies the specification property if for any $\epsilon>0$, there exists an integer $M\ge 0$ such that for any
$x_1,\cdots,x_k\in X$ and  any $n_1,\cdots,n_k\in\mathbb{N}$, there exists $y\in X$ such that
$d(f^i (x_j),f^{i+\sum_{t=1}^{j-1}(n_t+M)}(y))\le \epsilon$ holds for $0\le i\le n_j-1$ and $1\le j\le k$.
Here we set $\sum_{t=1}^0(n_t+M):=0$.
Roughly speaking, the specification property guarantees that the existence of an orbit, which traces all specified orbit-segments.
The one-way specification property, defined below, guarantees that the existence of an orbit, which traces specified orbit-segments ``contained in
the non-wandering set".
A point $x\in X$ is said to be non-wandering if for every neighborhood $U$ of $x$, there exists an integer $n>0$ such that $U\cap f^{-n}(U)\not=\emptyset$.
The set of all non-wandering points is called the non-wandering set and denoted by $\Omega(f)$.
It is well-known that $\mu(\Omega(f))=1$ for any $\mu\in\mathcal{M}_f(X)$ (see \cite[Proposition 6.19]{DGS} for instance).
Then we define the main property of this paper.

\begin{defi}
\label{low-sp}
We say that $(X,f)$ satisfies the one-way specification property if there exist compact subsets $X_1,\cdots,X_q$ of $X$ such that
\begin{enumerate}
\item $f(X_i)\subset X_i$ for $1\le i\le q$,

\item $\Omega(f)=\bigcup_{i=1}^qX_i$,

\item for any $\epsilon>0$, there exists an integer $M\ge 0$ such that for any $x_1\in X_{i(1)},\cdots,x_k\in X_{i(k)}$
with $i(1)\le\cdots\le i(k)$ and  any $n_1,\cdots,n_k\in\mathbb{N}$, there exists $y\in X$ such that
$$d(f^i (x_j),f^{i+\sum_{t=1}^{j-1}(n_t+M)}(y))\le \epsilon\hspace{1cm}(0\le i\le n_j-1,\ 1\le j\le k).$$
Here we set $\sum_{i=1}^0(n_t+M):=0$.
\end{enumerate}
\end{defi}

Our main examples appeared in \S3.2 do not satisfy the one-way specification property in general.
But they always satisfy the one-way (W)-specification property defined below, which is slightly weaker than the one-way specification property.

\begin{defi}
\label{low-wsp}
We say that $(X,f)$ satisfies the one-way (W)-specification property if there exist compact subsets $X_1,\cdots,X_q$ of $X$ such that
\begin{enumerate}
\item $f(X_i)\subset X_i$ for $1\le i\le q$,

\item $\Omega(f)=\bigcup_{i=1}^qX_i$,

\item for any $\epsilon>0$, there exists an integer $M\ge 0$ such that for any $x_1\in X_{i(1)},\cdots,x_k\in X_{i(k)}\in X$
with $i(1)\le\cdots\le i(k)$ and  any $n_1,\cdots,n_k\in\mathbb{N}$, there exist $y\in X$ and integers $M_1,\cdots,M_{k-1}\le M$ such that
$$d(f^i (x_j),f^{i+\sum_{t=1}^{j-1}(n_t+M_t)}(y))\le \epsilon\hspace{1cm}(0\le i\le n_j-1,\ 1\le j\le k).$$
Here we set $\sum_{i=1}^0(n_t+M_t):=0$.
\end{enumerate}
\end{defi}

The sense of ``one-way" in the above definitions seems to be well-understood once we give examples in \S3 (in particular,
see Figures \ref{fig:2}, \ref{fig:3} and \ref{fig:5}).
We sometimes say that $(X,f)$ satisfies the one-way specification property (or (W)-specification property) ``with $\{X_1,\cdots,X_q\}$".
We note that if $(X,f)$ satisfies the one-way (W)-specification property with $\{X_1,\cdots,X_q\}$,
then every ergodic measure is supported on some $X_i$ ($1\le i\le q$).
In fact, the following lemma holds.

\begin{lem}
\label{non-wandering}
Let $\mu\in\mathcal{M}_f^e(X)$ and $X_1,\cdots,X_q$ be compact subsets of $X$. Suppose that
$f(X_i)\subset X_i$ for $1\le i\le q$ and $\Omega(f)=\bigcup_{i=1}^qX_i$.
Then  there exists an integer $1\le i\le q$ such that $\mu(X_i)=1$.
\begin{proof}
Since $\mu(\Omega(f))=1$, there exists an integer $1\le i\le q$ such that $\mu(X_i)>0$.
Let $G_{\mu}$ be the set of all generic points for $\mu$, that is
$$G_{\mu}=\left\{x\in X:\lim_{n\rightarrow\infty}\mathcal{E}_n(x)=\mu\right\}$$
(see \cite[p.21]{DGS}). Since $\mu$ is ergodic, we have $\mu(G_{\mu})=1$, and which implies that
$\mu(G_{\mu}\cap X_i)>0$.
In particular, $G_{\mu}\cap X_i\not=\emptyset$.

Let us take a point $x\in G_{\mu}\cap X_i$. Then it follows from $f(X_i)\subset X_i$ that
$f^n(x)\in X_i$ for $n\ge 0$, which implies $(\mathcal{E}_n(x))(X_i)=1$ for any $n\ge 1$.
Since $\lim_{n\rightarrow\infty}\mathcal{E}_n(x)=\mu$ and $X_i$ is compact, we have
$$\mu(X_i)\ge\limsup_{n\rightarrow\infty}(\mathcal{E}_n(x))(X_i)=1,$$
which proves the lemma.
\end{proof}
\end{lem}

Let $\epsilon>0$ and $n\ge1$. A subset $E\subset X$ is called $(n,\epsilon)$-separated if for any two distinct points $x,y\in E$, there exists $0\le j\le n-1$ such that $d(f^j(x),f^j(y))>\epsilon$.
To show Theorem \ref{main1} , we use the following result, which is derived by a tiny modification of \cite[Proposition 2.1]{PfS}.

\begin{prop}
\label{PfSprop}
For any $\mu\in\mathcal{M}_f^e(X)$ and  any $h<h(\mu)$, there exists an $\epsilon>0$ such that for any neighborhood $F$ of $\mu$,
there exists an integer $N>0$ such that for any $n\ge N$ and for any Borel subset $B\subset X$ with
$\mu(B)=1$, there exists an $(n,\epsilon)$-separated set $\Gamma\subset \mathcal{E}_n^{-1}(F)\cap B$
such that $\sharp \Gamma\ge e^{nh}$.
\end{prop}

\subsection{Symbolic dynamics}
Let $A$ be a finite set and $A^{\mathbb{Z}_+}$ ($\mathbb{Z}:=\mathbb{N}\cup\{0\}$) be the set of all one-sided infinite sequences on the alphabet $A$,
endowed with the standard metric $d(\omega,\omega')=2^{-t(\omega,\omega')}$, where $t(\omega,\omega')=\min\{k:\omega_k\not=\omega'_k\}$.
The shift map $\sigma\colon A^{\mathbb{Z}_+}\to A^{\mathbb{Z}_+}$ is defined by $(\sigma(\omega))_i=\omega_{i+1}$ for $i\ge 0$.
We say that a compact subset $X\subset A^{\mathbb{Z}_+}$ is a subshift if $\sigma (X)\subset X$ holds.
We often denote $\sigma$ instead of the restriction map $\sigma|_X$ if no confusion arises.
The language of $X$, denoted by $\mathcal{L}=\mathcal{L}(X)$, is the set of all finite words that
appear in any sequence $\omega\in X$, that is,
$$\mathcal{L}(X):=\{w\in A^{\ast}:[w]\not=\emptyset\},$$
where $A^{\ast}=\bigcup_{n\ge 0}A^n$ and $[w]=\{\omega\in X:\omega_0\cdots \omega_{n-1}=w\}$ is a cylinder set for $w\in A^n$.
Given $w\in\mathcal{L}(X)$, let $|w|$ denote the length of $w$, and set $\mathcal{L}_n(X)=\{w\in\mathcal{L}(X):|w|=n\}$ for $n\ge 1$.
By our choice of the metric $d$ on $A^{\mathbb{Z}_+}$, we have the following proposition.

\begin{prop}
Let $(X,\sigma)$ be a subshift on a finite alphabet.
\begin{enumerate}
\item
$(X,\sigma)$ satisfies the one-way specification property if and only if
there exist subshifts $X_1,\cdots,X_q\subset X$ and an integer $M\ge 0$ such that $\Omega(\sigma|_X)=\bigcup_{i=1}^q\Omega(\sigma|_{X_i})$,
and for any $w_1\in\mathcal{L}(X_{i(1)}),\cdots,w_k\in\mathcal{L}(X_{i(k)})$ with $i(1)\le \cdots\le i(k)$, there exist $v_1,\cdots,v_{k-1}\in\mathcal{L}_M(X)$
such that $w_1v_1w_2v_2\cdots w_{k-1}v_{k-1}w_k\in\mathcal{L}(X)$.

\item
$(X,\sigma)$ satisfies the one-way (W)-specification property if and only if
there exist subshifts $X_1,\cdots,X_q\subset X$ and an integer $M\ge 0$ such that $\Omega(\sigma|_X)=\bigcup_{i=1}^q\Omega(\sigma|_{X_i})$,
and for any $w_1\in\mathcal{L}(X_{i(1)}),\cdots,w_k\in\mathcal{L}(X_{i(k)})$ with $i(1)\le \cdots\le i(k)$, there exist $v_1,\cdots,v_{k-1}\in\mathcal{L}(X)$
such that $|v_i|\le M$ for $1\le i\le k-1$ and $w_1v_1w_2v_2\cdots w_{k-1}v_{k-1}w_k\in\mathcal{L}(X)$.
\end{enumerate}
\end{prop}

\subsection{Symbolic spaces of piecewise expanding interval systems}
Let $I=[0,1]$ be the unit interval, and let $f\colon I\to I$
be such that there exist $p\ge 2$ and $0=a_0<a_1<\cdots<a_p=1$
such that writing $I_j=(a_j,a_{j+1})$, the restriction
$f|_{I_j}$ is $C^1$ and satisfies $\alpha\le |f'(x)|\le \beta$
$(x\in\bigcup_{j=1}^{p-1}I_j)$ for some $1<\alpha\le\beta$.
We say that $(I,f)$ is a piecewise expanding interval system.

Let $A=\{0,\cdots,p-1\},\ S=\{a_0,\cdots,a_p\}$ and $I'=I\setminus\bigcup_{i\ge 0}
f^{-i}S$.
We define the map $i'\colon I'\to A^{\mathbb{Z}_+}$
by $(i'(x))_k=j$ if $f^k(x)\in I_j$.
We note that $i'$ is well-defined and injective on $I'$
since $|f'(x)|\ge \alpha>1$ for all $x\in I\setminus S$.
We set $\Sigma_f={\rm cl}\{i'(I')\}$, where ${\rm cl}(A)$
denotes the closure of the set $A$.
Then $(\Sigma_f,\sigma)$ is a subshift.
We call $(\Sigma_f,\sigma)$ the symbolic space of $(I,f)$.

We define $\Psi\colon\Sigma_f\to I$ by
$\Psi (\omega)=\bigcap_{k=0}^{\infty}f^{-k}{\rm cl}(I_{\omega_k})$.
It is well-known that $(I,f)$ is a topological
factor of $(\Sigma_f,\sigma)$ with $\Psi$ as a factor map.
%Moreover, there exists a Borel subset $Y\subset \Sigma_f$ such that both $\Sigma_f\setminus Y$ and
%$[0,1]\setminus\Psi(Y)$ are countable and $\Psi\colon Y\to\Psi(Y)$ is bijective.

\section{Examples}
\subsection{Simple examples}
In this subsection, we give simple examples of symbolic and non-symbolic systems, which satisfies the
one-way specification property.
\begin{ex}
\label{se1}
We consider a piecewise expanding interval system whose symbolic space has the one-way specification property.
Let $I=[0,1]$ be the unit interval and define $f\colon I\to I$ by
$$f(x)=
\left\{
\begin{array}{ll}
3x & (0\le x<\frac{1}{6}); \vspace{0.1cm}\\
3x-\frac{1}{2} & (\frac{1}{6}\le x<\frac{1}{2});\vspace{0.1cm}\\
3x-1 & (\frac{1}{2}\le x<\frac{2}{3});\vspace{0.1cm}\\
3x-\frac{3}{2} & (\frac{2}{3}\le x<\frac{5}{6});\vspace{0.1cm}\\
3x-2 & (\frac{5}{6}\le x\le 1).
\end{array}
\right.
$$
(A graph of $f$ is sketched in Figure \ref{fig:1}.)
\begin{figure}[htbp]
 \begin{center}
  \includegraphics[width=65mm]{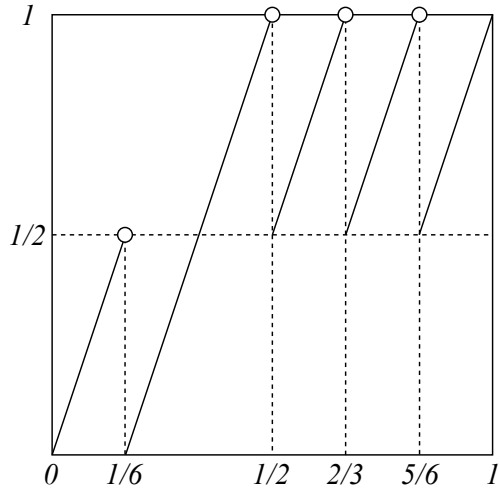}
 \end{center}
 \caption{A graph of $f$.}
 \label{fig:1}
\end{figure}
Let $(\Sigma_f,\sigma)$ be the symbolic space of $(I,f)$.
Then it is not difficult to see that the language of $\Sigma_f$ is the set of all
subwords of the form $w1v$ $\displaystyle(w\in\{0,1\}^{\mathbb{Z}_+},\ v\in\{2,3,4\}^{\mathbb{Z}_+})$.
In other words, $\Sigma_f$ is the set of all infinite paths started from some vertex of a Graph $\Gamma$, which is shown in Figure \ref{fig:2}.
\begin{figure}[htbp]
 \begin{center}
  \includegraphics[width=80mm]{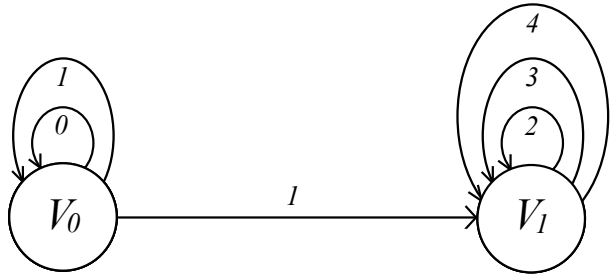}
 \end{center}
 \caption{A graph $\Gamma$.}
 \label{fig:2}
\end{figure}
Thus, $(\Sigma_f,\sigma)$ satisfies the one-way specification property with $\{\{0,1\}^{\mathbb{Z}_+},\{2,3,4\}^{\mathbb{Z}_+}\}$.
\end{ex}

\begin{ex}
\label{se2}
Let $S^1$ be the circle $\mathbb{R}/\mathbb{Z}$ and take a system of representatives $[0,1)$ of $S^1$.
We define $f\colon S^1\to S^1$ by
$$f(\theta)=\theta+\frac{1}{10}\sin (2\pi\theta).$$
In Figure \ref{fig:3}, we sketch a phase portrait of $f$.
\begin{figure}[htbp]
 \begin{center}
  \includegraphics[width=65mm]{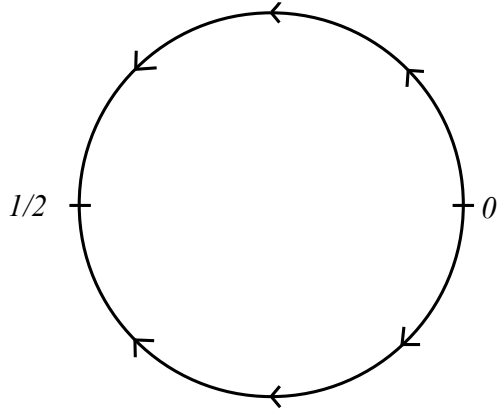}
 \end{center}
 \caption{A phase portrait of $f$. The non-wandering set of $f$ consists of only two points $0$ and $\frac{1}{2}$.
$0$ is a unique source and $\frac{1}{2}$ is a unique sink.}
 \label{fig:3}
\end{figure}
Then it is known that $\Omega(f)=\{0,\frac{1}{2}\}$ (see \cite[Example 12.1 in \S8.12]{Ro} for instance) and it is not difficult to see that $f$ satisfies the one-way specification property
with $\{\{0\},\{\frac{1}{2}\}\}$.
\end{ex}

\subsection{($-\beta$)-transformations}
For $\beta>1$, we define a $(-\beta)$-transformation $T_{-\beta}\colon [0,1]\to [0,1]$ as follows.
Let $b:=\max\{k\in\mathbb{Z}:k<\beta\}$ and set $I_i:=\left(\frac{i}{\beta},\frac{i+1}{\beta}\right)$ for $0\le i\le b-1$
and $I_b:=\left(\frac{b}{\beta},1\right)$.
First, we define values of $T_{-\beta}$ on $\bigcup_{i=0}^bI_b$ by
$$T_{-\beta}(x):=-\beta x+(i+1)\hspace{1cm}(0\le i\le b).$$
Then, for $x\not \in \bigcup_{i=0}^{\infty}T_{-\beta}^{-i}\left(\left\{0,\frac{1}{\beta},\cdots,\frac{b}{\beta},1\right\}\right)$, one can define
$i'(x)\in\{0,1,\cdots,b\}^{\mathbb{Z}_+}$ by $(i'(x))_k=j$ if $T_{-\beta}^k(x)\in I_j$.
We set $\displaystyle i'(1):=\lim_{x\rightarrow 1-0}i'(x)$.

In what follows we will define values of $T_{-\beta}$ on $\left\{0,\frac{1}{\beta},\cdots\frac{b}{\beta},1\right\}$.
\vspace{0.1cm}\\
{\em Case} 1. $i'(1)$ is of the form $(w0)^{\infty}=w0w0w0\cdots$ for some $w\in\{0,1,\cdots,b\}^{\ast}$.
\vspace{0.1cm}\\
In this case, we define $T_{-\beta}$ on $\left\{0,\frac{1}{\beta},\cdots\frac{b}{\beta},1\right\}$ by
$$\displaystyle T_{-\beta}(x)=
\left\{
\begin{array}{ll}
1 & (x=0); \\
0 & (x=\frac{1}{\beta},\cdots,\frac{b}{\beta}); \\
\lim_{x\rightarrow 1-0}T_{-\beta}(x) & (x=1).
\end{array}
\right.
$$

\begin{figure}[htbp]
 \begin{center}
  \includegraphics[width=120mm]{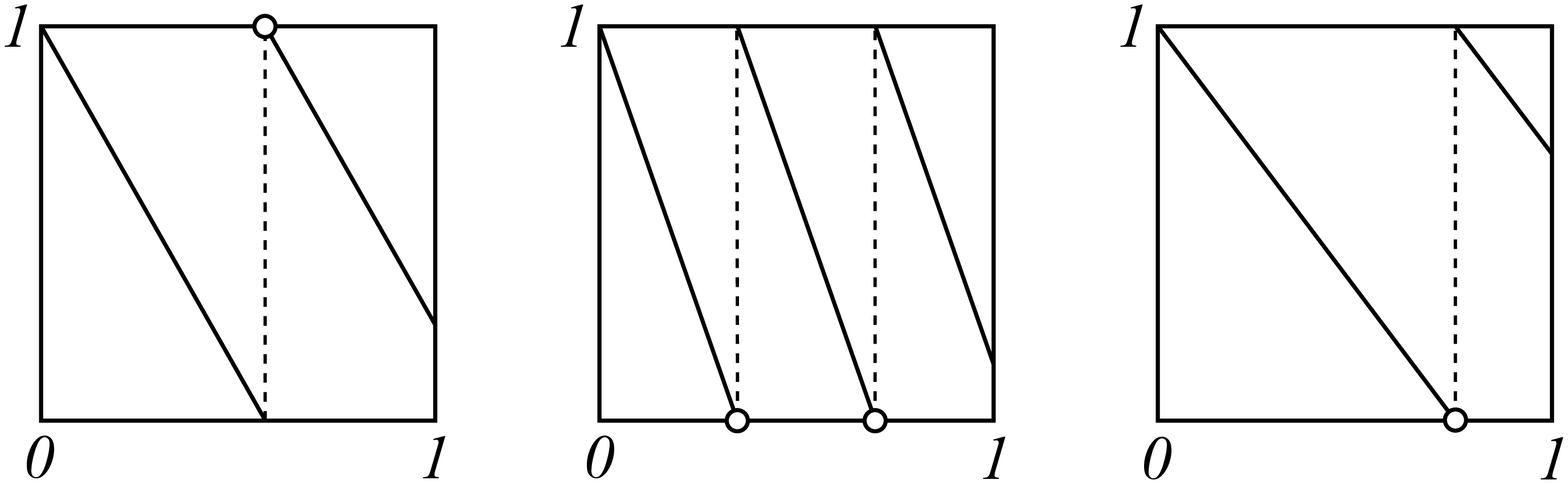}
 \end{center}
 \caption{Graphs of $T_{-\beta}$ for $i'(1)=(1000)^{\infty}$ (left), $i'(1)=(2011)^{\infty}$ (center) and $i'(1)=100(1)^{\infty}$ (right).}
 \label{fig:4}
\end{figure}

\noindent
{\em Case} 2. $i'(1)$ is not of the form $(w0)^{\infty}=w0w0w0\cdots$.
\vspace{0.1cm}\\
In this case, we define $T_{-\beta}$ on $\left\{0,\frac{1}{\beta},\cdots\frac{b}{\beta},1\right\}$ by
$$\displaystyle T_{-\beta}(x)=
\left\{
\begin{array}{ll}
1 & (x=0); \\
1 & (x=\frac{1}{\beta},\cdots,\frac{b}{\beta}); \\
\lim_{x\rightarrow 1-0}T_{-\beta}(x) & (x=1).
\end{array}
\right.
$$
(The situations are sketched in Figure \ref{fig:4} for $i'(1)=(1000)^{\infty}$, $i'(1)=(2011)^{\infty}$
and $i'(1)=100(1)^{\infty}$.)

We define the alternating order $<_{{\rm alt}}$ on $\{0,1,\cdots,b\}^{\mathbb{Z}_+}$ by
$$(\omega_i)<_{{\rm alt}}(\omega'_i)\text{ if and only if }
x_i=y_i\ (0\le i\le n-1)\text{ and }
\left\{
\begin{array}{ll}
\omega_n<\omega'_n & (n\text{ is odd}); \\
\omega_n>\omega'_n & (n\text{ is even})
\end{array}
\right.$$
hold. We denote $(\omega_i)\le_{{\rm alt}}(\omega'_i)$ if either $(\omega_i)<_{{\rm alt}}(\omega'_i)$ or $(\omega_i)=(\omega'_i)$ holds.
Let $(\Sigma_{-\beta},\sigma)$ be the symbolic space of $([0,1],T_{-\beta})$.
Then following the proofs appeared in \cite[\S2]{IS} with a tiny modification, we have
\begin{equation}
\label{eq1}
\Sigma_{-\beta}=\{(\omega_i):\sigma^n((\omega_i))\le_{{\rm alt}}i'(1)\text{ for }n\ge 0\}.
\end{equation}
We denote $i'(1)=s_0s_1s_2\cdots$ and let $\Gamma_{-\beta}$ be a graph constructed by the following rule.
\begin{itemize}
\item
there is an edge $V_i\xrightarrow{s_i} V_{i+1}$ for any $i\ge 0$.

\item
If $i$ is even, $0\le a\le s_i-1$, and $s_0\cdots s_{j-1}$ is the suffix of maximal length of $s_0\cdots s_{i-1}a$,
then there is an edge $V_i\xrightarrow{a}V_j$.

\item
If $i$ is odd, $s_i+1\le b\le s_0$, and $s_0\cdots s_{j-1}$ is the suffix of maximal length of $s_0\cdots s_{i-1}b$,
then there is an edge $V_i\xrightarrow{b}V_j$.

\item
If $i$ is odd and $s_i<s_0$, then there is an edge $V_i\xrightarrow{s_0} V_1$.
\end{itemize}

\begin{figure}[htbp]
 \begin{center}
  \includegraphics[width=125mm]{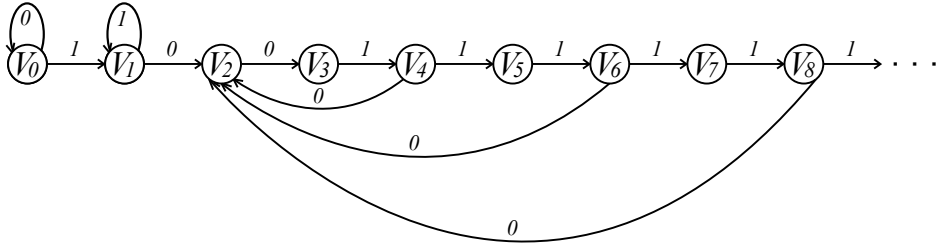}
 \end{center}
 \caption{A graph $\Gamma_{-\beta}$ for $i'(1)=100(1)^{\infty}$.}
 \label{fig:5}
\end{figure}

The situation is sketched in Figure \ref{fig:5} for $i'(1)=100(1)^{\infty}$.
Then by the equation (\ref{eq1}) and the definition of the alternating order, we have
$$\Sigma_{-\beta}=\{\text{Set of all one-sided infinite paths started from some vertex }V_i\text{ of }\Gamma_{-\beta}\}.$$

In what follows we will show the following proposition.
\begin{prop}
\label{beta-OWSP}
Suppose that $i'(1)$ is eventually periodic. Then $(\Sigma_{-\beta},\sigma)$ satisfies the one-way (W)-specification property.
\begin{proof}
We denote
$i'(1)=s_0s_1s_2\cdots=t_0\cdots t_{u-1}(t_u\cdots t_{u+v-1})^{\infty}$ for some $u\ge 0$ and $v\ge 1$.
Since $i'(1)=\lim_{x\rightarrow 1-0}i'(x)=s_0s_1\cdots$, the cylinder set $[s_0\cdots s_i]$ is not a single point set for any $i\ge 0$.
Thus, there is an edge
$V_i\xrightarrow{a}V_j$ for some $i\ge u+v+2$, $j<i$ and $0\le a\le s_0$.
\begin{lem}
\label{beta-lem1}
We have $j\le u+v$. Moreover, for each $n\ge 0$, there is an edge $V_{i+2nv}\xrightarrow{a}V_j$.
\begin{proof}
First we show $j\le u+v$. By contradiction, suppose that $j\ge u+v+1$. 
Since there is an edge $V_i\xrightarrow{a}V_j$, $s_0\cdots s_{j-1}$ is the suffix of maximal length of $s_0\cdots s_{i-1}a$.
Thus we can find a word $w'$ such that
$s_0\cdots s_{i-1}a=w's_0\cdots s_{j-1}$. Since $j\ge u+v+1$, we have $s_{j-1}=s_{j-v-1}$.
Therefore we have $a=s_{j-1}=s_{j-v-1}=s_{i-v}=s_i$, which contradicts with $s_i\not= a$.

Let $n\ge 0$. Since $s_0\cdots s_{j-1}$ is the suffix of $s_0\cdots s_{i-1}a$ and $j\le u+v$,
$s_0\cdots s_{j-1}$ is also the suffix of $s_0\cdots s_{i+2nv-1}a$.
In what follows we will show that $s_0\cdots s_{j-1}$ is the suffix of maximal length of $s_0\cdots s_{i+2nv-1}a$.
By contradiction, we assume that there exists a number $j'>j$ such that
$s_0\cdots s_{j'-1}$ is the suffix of $s_0\cdots s_{i+2nv-1}a$.
Then $s_0\cdots s_j$ is the suffix of $s_0\cdots s_{i-1}a$.
However, $s_0\cdots s_{j-1}$ is the suffix of maximal length of $s_0\cdots s_{i-1}a$.
This is a contradiction.
\end{proof}
\end{lem}
For integers $0\le n\le m$, let $\Gamma_{n,m}$ (resp. $\Gamma_{n,\infty}$) be a subgraph of $\Gamma_{-\beta}$
from $V_n$ to $V_m$ (resp. from $V_n$), that is,  the vertexes of $\Gamma_{n.m}$ (resp. $\Gamma_{n,\infty}$) are
$V_n,V_{n+1},\cdots,V_m$ (resp. $V_n,V_{n+1},\cdots$) and the edges of $\Gamma_{n,m}$ are the set of all edges of $\Gamma_{-\beta}$
from $V_i$ to $V_j$ for some $n\le i,j\le m$ (resp. $i,j\ge n$). We set
$$N:=\min\{k\ge 0:\Gamma_{k,\infty}\text{ is irreducible}\}.$$
By Lemma \ref{beta-lem1}, such $N$ exists. We define sequences of integers $\{l_i\}$ and $\{n_i\}$ inductively as follows.
First, we set $l_1:=0$ and
$$n_1:=\max\{l_1\le k< N:\Gamma_{l_1,k}\text{ is irreducible}\}.$$
Next, we set
$$l_2:=\min\{k> n_1:\text{there are at least two edges into $V_k$}\}\text{ and}$$
$$n_2:=\{l_2\le k< N:\Gamma_{l_2,k}\text{ is irreducible}\}.$$
Using this procedure inductively, we can construct integers $0=l_1\le n_1<l_2\le n_2<\cdots<l_{q-1}\le n_{q-1}<l_q=N<n_q=\infty$ such that
\begin{itemize}
\item
$\Gamma_{l_i,n_i}$ is irreducible for each $1\le i\le q$;

\item
$V_k\xrightarrow{s_k}V_{k+1}$ is only edge into $V_k$ for $n_i\le k \le l_{i+1}$ and $1\le i\le q-1$;

\item
there is no edge from $V_k$ ($k\not \in \{l_i,\cdots,n_i\}$) to $\Gamma_{l_i,n_i}$ except
$V_{l_i-1}\xrightarrow{s_{l_i-1}}V_{l_i}$ for $1\le i\le q$.
\end{itemize}

(In Figure \ref{fig:6}, we sketched the situation for $i'(1)=100(1)^{\infty}$.)
\begin{figure}[htbp]
 \begin{center}
  \includegraphics[width=125mm]{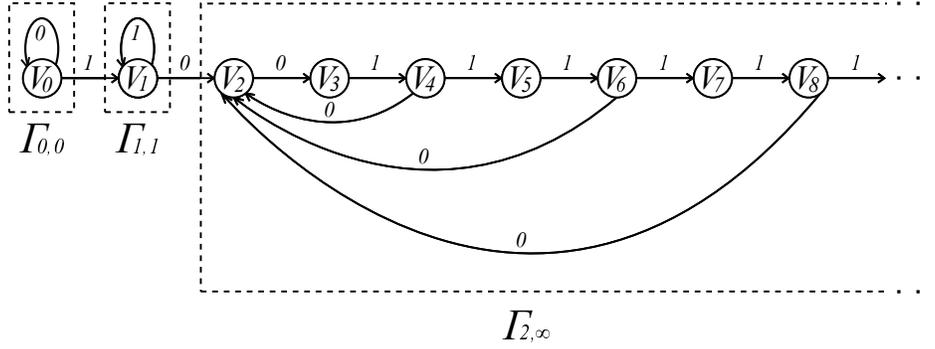}
 \end{center}
 \caption{A graph representation of $(\Sigma_{-\beta},\sigma)$ for $i'(1)=100(1)^{\infty}$.}
 \label{fig:6}
\end{figure}

Let $X_i$ be a subshift generated by the graph $\Gamma_{l_i,n_i}$ for $1\le i\le q$.
More precisely, $X_i$ is the set of all infinite paths started from some vertex in $\Gamma_{l_i,n_i}$.
In what follows we will show that $(\Sigma_{-\beta},\sigma)$ satisfies the one-way (W)-specification property
with $\{X_1,\cdots,X_q\}$.
Since $i'(1)$ is eventually periodic, we can show that $X_q$ is sofic in a similar way to the proof of \cite[Proposition 2]{FL}.
Thus, $X_i$ is sofic and irreducible for each $1\le i\le q$. Thus, $X_i$ satisfies the (W)-specification property,
that is, there exists an integer $M\ge 0$ such that for any $w_1,\cdots,w_k\in\mathcal{L}(X_i)$, we can find
$v_1,\cdots,v_{k-1}\in\mathcal{L}(X_i)$ such that $|v_j|\le M$ for $1\le j\le k-1$ and
$$w_1v_1w_2v_2\cdots w_{k-1}v_{k-1}w_k\in\mathcal{L}(X_i).$$
Therefore, $(\Sigma_{-\beta},\sigma)$ satisfies the one-way (W)-specification property with $\{X_1,\cdots,X_q\}$.
\end{proof}
\end{prop}

\begin{rem}
\label{rem-minus-beta}
As we mentioned in \S1, our definition of the $(-\beta)$-transformation is slightly different from the original one in \cite{IS}.
Ito and Sadahiro introduced a $(-\beta)$-transformation $T_{-\beta}^{{\rm IS}}$ on the interval $[\frac{-\beta}{\beta+1},\frac{1}{\beta+1})$ defined by
$T_{-\beta}^{{\rm IS}}(x):=-\beta x-\lfloor \frac{\beta}{\beta+1}-\beta x\rfloor,$
where $\lfloor x\rfloor$ denotes the largest integer not exceeding a real number $x$.
After that, Liao and Steiner \cite{LS} considered a $(-\beta)$-transformation $T_{-\beta}^{{\rm LS}}$ on the interval $(0,1]$ defined by
$T_{-\beta}^{{\rm LS}}(x):=-\beta x+\lfloor \beta x\rfloor +1,$
which is topologically conjugate to $T_{-\beta}^{{\rm IS}}$. Then, the differences between our $T_{-\beta}$ and Liao-Steiner's $T_{-\beta}^{{\rm LS}}$
are the following.
\begin{itemize}
\item
$T_{-\beta}$ is defined on $[0,1]$, but $T_{-\beta}^{{\rm LS}}$ is defined on $(0,1]$.

\item
In {\em Case} 1, that is, $i'(1)$ is of the form $(w 0)^{\infty}$, then $T_{-\beta}(x)=0$ but $T_{-\beta}^{{\rm LS}}(x)=1$ for $x=\frac{1}{\beta},\cdots,\frac{b}{\beta}$.
If we further assume $\beta>1$ is  an integer, then $T_{-\beta}(1)=0$ but $T_{-\beta}^{{\rm LS}}(1)=1$.
\end{itemize}
The reason why we give a slightly modified version of the definition of the ($-\beta$)-transformation is the following.
\begin{itemize}
\item
To formulate the large deviation principle, one usually requires that the spece is Polish.

\item
By our definition of $T_{-\beta}$, $\Psi([w])$ is an interval with positive length for any $w\in\mathcal{L}(\Sigma_{-\beta})$, where
$\Psi$ is a map defined in \S2.3. If we define $T_{-\beta}$ as a natural extension of $T_{-\beta}^{{\rm LS}}$,
then it may occur that $\Psi([w])$ is a single point for some $w\in\Sigma_{-\beta}$.
\end{itemize}
We note that in {\em Case} 2, $T_{-\beta}(x)=-\beta x+\lfloor \beta x\rfloor +1$ and so our definition of the $(-\beta)$-transformation $T_{-\beta}$
is a natural extension of $T_{-\beta}^{{\rm LS}}$. Therefore,
the symbolic space $(\Sigma_{-\beta},\sigma)$ of $([0,1],T_{-\beta})$ coincides with that of $((0,1],T_{-\beta}^{{\rm LS}})$ (and so, also with that of
$([\frac{-\beta}{\beta+1},\frac{1}{\beta+1}),T_{-\beta}^{{\rm IS}})$).
\end{rem}

\section{Proof of Theorem \ref{main1}}
In this section, we give a proof of Theorem \ref{main1}. First we establish the upper bound.
Since we assume \ref{A.III} and \ref{A.IIII},
it follows from \cite[Theorems 3.1 and 3.2]{PfS} that
$$\displaystyle\limsup_{n\rightarrow\infty}\frac{1}{n}\log m(\{x\in X:\mathcal{E}_n(x)\in F\})\le\sup_{\mu\in F}\mathcal{F}_{\varphi}(\mu)$$
holds for any closed set $F\subset\mathcal{M}(X)$.

To get the lower bound, it is sufficient to show that for any $\mu\in\mathcal{M}_f(X)$
and any open neighborhood $G\subset\mathcal{M}(X)$ of $\mu$,
\begin{equation}
\label{lb}
\displaystyle\liminf_{n\rightarrow\infty}\frac{1}{n}\log
m(\{x\in X:\mathcal{E}_n(x)\in G\})\ge h(\mu)-\int \varphi d\mu
\end{equation}
holds.
Let $\mu\in\mathcal{M}_f(X)$, $G$ be an open neighborhood of $\mu$ and $\eta>0$.
Choose $\zeta>0$ so small that $D(\mu,\nu)\le 5\zeta$ implies that $\nu\in G$ and $|\int\varphi d\mu-\int\varphi d\nu|\le\eta$.
It follows from \ref{A.II} that there exist $\epsilon_1>0$ and $N_1\in\mathbb{N}$
such that for any $0<\epsilon\le \epsilon_1$, for any $n\ge N_1$ and for any $x\in X$,
$$m(B_n(x,\epsilon))\ge \exp\left(-\sum_{j=0}^{n-1}\varphi(f^j(x))-n\eta\right)$$
holds.
We can also find $0<\epsilon_2\le\epsilon_1$ such that $D(\delta_x,\delta_y)\le\zeta$ whenever $d(x,y)\le2\epsilon_2$.

On the other hand, by Ergodic Decomposition Theorem and the affinity of the
entropy map, there exists $\lambda=\sum_{i=1}^pa_i\mu_i$ such that
\begin{itemize}
\item
the $\mu_i$ are ergodic;

\item
the $a_i$ are real numbers in $(0,1]$ such that $\sum_{i=1}^pa_i=1$;

\item
$D(\mu,\lambda)\le \zeta$;

\item
$h(\lambda)>h(\mu)-\eta$.
\end{itemize}

By the condition \ref{A.I}, $(X,f)$ satisfies the one-way (W)-specification property with some $\{X_1,\cdots,X_q\}$.
Thus, it follows from Lemma \ref{non-wandering} that for each $1\le i\le p$, we can find
$1\le j(i)\le q$ such that $\mu_i(X_{j(i)})=1$.
Without loss of generality, we may assume $j(1)\le j(2)\le \cdots\le j(p)$.
By Proposition \ref{PfSprop}, we can find $0<\epsilon\le \epsilon_2$ and $N_2\ge N_1$ such that
for any $n\ge N_2$ and any $1\le i\le p$, there exists an $(n_i,3\epsilon)$-separated subset
$\Gamma_{n,i}\subset\mathcal{E}_{n_i}^{-1}(\mathbb{B}(\mu_i,\zeta))\cap X_{j(i)}$ such that
$$\sharp\Gamma_{n,i}\ge \exp(n_i(h(\mu_i)-\eta))$$
holds. Here we set $\mathbb{B}(\mu_i,\eta):=\{\nu\in\mathcal{M}(X):D(\mu_i,\nu)\le\eta\}$ and $n_i:=\max\{k\in\mathbb{Z}:a_in\ge k\}$.

Let $M$ be as in the definition of the one-way (W)-specification property. Choose an integer $N\ge N_2$
so large that $\frac{1}{N}\log(M+1)^{p-1}\le \eta$ and
$$\sum_{i=1}^p\left|a_i-\frac{n_i}{n_1+\cdots+n_p}\right|\le\zeta$$
hold.
Fix any integer $n\ge N$. Then by the one-way (W)-specification property, for any
$\mathbb{X}=(x_1,\cdots,x_p)\in\prod_{i=1}^p\Gamma_{n,i}$, we can find $y(\mathbb{X})\in X$ and
integers $M_1(\mathbb{X}),\cdots,M_{p-1}(\mathbb{X})\le M$ such that
$$d(f^j(x_i),f^{j+\sum_{t=1}^{i-1}(n_t+M_t(\mathbb{X}))}(y(\mathbb{X})))\le \epsilon\hspace{1cm}(0\le j\le n_i-1,\ 1\le i\le p).$$
Since $0\le M_i(\mathbb{X})\le M$ for $1\le i\le p$, it is easy to see that there exist integers
$0\le M_1,\cdots,M_{p-1}\le M$ such that
$$\sharp\left\{\mathbb{X}\in\prod_{i=1}^p\Gamma_{n,i}:M_i(\mathbb{X})=M_i,\ 1\le i\le p-1\right\}\ge \frac{1}{(M+1)^{p-1}}\sharp\prod_{i=1}^p\Gamma_{n,i}.$$
For the notational simplicity, we set
$$\Gamma_n:=\left\{\mathbb{X}\in\prod_{i=1}^p\Gamma_{n,i}:M_i(\mathbb{X})=M_i,\ 1\le i\le p-1\right\}.$$
Then we have
\begin{eqnarray*}
\sharp\Gamma_n
&\ge&
 \frac{1}{(M+1)^{p-1}}\sharp\prod_{i=1}^p\Gamma_{n,i} \\
&\ge&
\frac{1}{(M+1)^{p-1}}\prod_{i=1}^p\exp(n_i(h(\mu_i)-\eta)) \\
&\ge&
\frac{1}{(M+1)^{p-1}}\exp\left(\sum_{i=1}^pa_in(h(\mu_i)-2\eta)\right) \\
&\ge&
\frac{1}{(M+1)^{p-1}}\exp(n(h(\mu)-3\eta)) \\
&\ge&
\exp(n(h(\mu)-4\eta)).
\end{eqnarray*}
We set $\tilde{n}:=\sum_{i=1}^{p-1}(n_i+M_i)+n_p$. Taking $N$ large if necessary, we may assume that
$$D\left(\mathcal{E}_n(x),\sum_{i=1}^p\frac{n_i}{n_1+\cdots+n_p}\mathcal{E}_{n_i}(f^{\sum_{t=1}^{i-1}(n_t+M_t)}(x))\right)\le\zeta$$
for $x\in X$, $\frac{\tilde{n}-n}{n}\int\varphi d\mu\le\eta$, $\frac{2\tilde{n}}{n}\eta\le 3\eta$ and $\tilde{n}\ge \max\{N_1,N_2\}$ hold.

\begin{lem}
\label{mainlem1}
For any $\mathbb{X},\mathbb{X}'\in\Gamma_n$ with $\mathbb{X}\not=\mathbb{X}'$, we have
$$B_{\tilde{n}}(y(\mathbb{X}),\epsilon)\cap B_{\tilde{n}}(y(\mathbb{X}'),\epsilon)=\emptyset.$$
\begin{proof}
We denote $\mathbb{X}=(x_1,\cdots,x_p)$ and $\mathbb{X}'=(x_1',\cdots,x_p')$.
By $\mathbb{X}\not=\mathbb{X}'$, we can find $1\le i\le p$ such that $x_i\not=x_i'$.
Since $x_i,x_i'\in\Gamma_{n,i}$ and $\Gamma_{n,i}$ is $(n_i,3\epsilon)$-separated, there exists $0\le j\le n_i-1$ such that
$d(f^j(x_i),f^j(x_i'))>3\epsilon$. Thus, we have
\begin{eqnarray*}
{}
& &
d(f^{j+\sum_{t=1}^{i-1}(n_t+M_t)}(y(\mathbb{X})),f^{j+\sum_{t=1}^{i-1}(n_t+M_t)}(y(\mathbb{X}'))) \\
&\ge&
d(f^j (x_i),f^j(x_i'))-d(f^j(x_i),f^{j+\sum_{t=1}^{i-1}(n_t+M_t)}(y(\mathbb{X})))\\
& &
\hspace{2.8cm}-d(f^j(x_i'),f^{j+\sum_{t=1}^{i-1}(n_t+M_t)}(y(\mathbb{X}'))) \\
&>&
3\epsilon-2\epsilon \\
&=&
\epsilon,
\end{eqnarray*}
which proves the lemma.
\end{proof}
\end{lem}

\begin{lem}
\label{mainlem2}
Suppose that $\mathbb{X}\in\Gamma_n$ and $x\in B_{\tilde{n}}(y(\mathbb{X}),\epsilon)$. Then we have $D(\mu,\mathcal{E}_n(x))\le 5\zeta$.
In particular, $\mathcal{E}_n(x)\in G$ and $\left|\int\varphi d\mu-\frac{1}{n}\sum_{j=0}^{n-1}\varphi(f^j(x))\right|\le \eta$ hold.
\begin{proof}
Denote $\mathbb{X}=(x_1,\cdots,x_p)$ and suppose that $x\in B_{\tilde{n}}(y(\mathbb{X}),\epsilon)$. Then we have
$$d(f^j(x_i),f^{j+\sum_{t=1}^{i-1}(n_t+M_t)}(x))\le 2\epsilon\hspace{1cm}(0\le j\le n_i-1,\ 1\le i\le p).$$
Thus, it follows from $x_i\in\mathcal{E}_{n_i}^{-1}(\mathbb{B}(\mu_i,\zeta))$ for $1\le i\le p$ and Lemma \ref{weak metric} that
\begin{eqnarray*}
D(\mu,\mathcal{E}_n(x))
&\le&
D\left(\sum_{i=1}^pa_i\mu_i,\sum_{i=1}^p\frac{n_i}{n_1+\cdots+n_p}\mathcal{E}_{n_i}(f^{\sum_{t=1}^{i-1}(n_t+M_t)}(x))\right)+2\zeta \\
&\le&
D\left(\sum_{i=1}^pa_i\mu_i,\sum_{i=1}^p\frac{n_i}{n_1+\cdots+n_p}\mathcal{E}_{n_i}(x_i)\right)+3\zeta \\
&\le&
5\zeta,
\end{eqnarray*}
which proves the lemma.
\end{proof}
\end{lem}

Now we continue the proof of Theorem \ref{main1}. By Lemmas \ref{mainlem1} and \ref{mainlem2}, we have
\begin{eqnarray*}
{}
& &
\liminf_{n\rightarrow\infty}\frac{1}{n}\log m(\{x\in X:\mathcal{E}_n(x)\in G\}) \\
&\ge&
\liminf_{n\rightarrow\infty}\frac{1}{n}\log m\left(\bigcup_{\mathbb{X}\in\Gamma_n}B_{\tilde{n}}(y(\mathbb{X}),\epsilon)\right) \\
&=&
\liminf_{n\rightarrow\infty}\frac{1}{n}\log\sum_{\mathbb{X}\in\Gamma_n} m\left(B_{\tilde{n}}(y(\mathbb{X}),\epsilon)\right) \\
&\ge&
\liminf_{n\rightarrow\infty}\frac{1}{n}\log\sharp\Gamma_n\exp\left(-\sum_{j=0}^{\tilde{n}-1}\varphi(f^j(y(\mathbb{X})))-\tilde{n}\eta\right) \\
&\ge&
\liminf_{n\rightarrow\infty}\frac{1}{n}\exp\left(n(h(\mu)-4\eta)-\tilde{n}\int\varphi d\mu-2\tilde{n}\eta\right) \\
&\ge&
h(\mu)-\int\varphi d\mu -8\eta.
\end{eqnarray*}
This implies the equation (\ref{lb}).

\section{Proof of the corollaries}
\subsection{Proof of Corollary \ref{maincor}}
In this subsection, we give a proof of Corollary \ref{maincor}.
It is well-known that if $X$ is a subshift, then the entropy map $\mu\mapsto h(\mu)$ is upper-semicontinuous on $\mathcal{M}_{\sigma}(X)$.
Thus, Theorem \ref{main1} together with \cite[Proposition 4.2]{PfS} imply the corollary.

\subsection{Proof of Corollary \ref{p-e}}
The purpose of this subsection is to prove Corollary \ref{p-e}.
For $\tilde{\mu}\in\mathcal{M}(\Sigma_f)$, we set $\hat{\Psi}(\tilde{\mu}):=\tilde{\mu}\circ\Psi^{-1}$.
To prove Corollary \ref{p-e}, we use the following lemmas.

\begin{lem}
\label{conti}
The map $\hat{\Psi}\colon \mathcal{M}(\Sigma_f)\to \mathcal{M}([0,1])$ is continuous.
\begin{proof}
Suppose that $\{\tilde{\mu}_n\}\subset\mathcal{M}(\Sigma_f)$ and
$\displaystyle\lim_{n\rightarrow\infty}\tilde{\mu}_n=\tilde{\mu}$.
Then for any continuous function $\psi\colon [0,1]\to \mathbb{R}$, we have
\begin{eqnarray*}
\lim_{n\rightarrow\infty}\int\psi d\hat{\Psi}(\tilde{\mu}_n)
&=&
%\lim_{n\rightarrow\infty}\int\psi d\tilde{\mu}_n\circ\Psi^{-1} \\
%&=&
\lim_{n\rightarrow\infty}\int\psi\circ\Psi d\tilde{\mu}_n \\
&=&
\int\psi \circ\Psi d\tilde{\mu} \\
&=&
\int\psi d\hat{\Psi}(\tilde{\mu}),
\end{eqnarray*}
which implies that $\hat{\Psi}$ is continuous.
\end{proof}
\end{lem}

By Corollary \ref{maincor} and the condition \ref{2.2},
$(\Sigma_f,\sigma)$ satisfies a level-2 large deviation principle with $L_f$ and
the rate function is the free energy $\widetilde{\mathcal{F}}\colon\mathcal{M}(\Sigma_f)\to[-\infty,0]$ with respect to $\varphi\circ \Psi$, that is,
$$\widetilde{\mathcal{F}}(\mu)=
\left\{
\begin{array}{ll}
h(\mu)-\int\varphi\circ\Psi d\mu & (\mu\in\mathcal{M}_{\sigma}(\Sigma_f)); \\
-\infty & (\text{otherwise}).
\end{array}
\right.
$$
Let $\mathcal{F}_{\varphi}\colon\mathcal{M}([0,1])\to [-\infty,0]$ be the free energy with respect to $\varphi$ given by (\ref{free}).
\begin{lem}
\label{com}
For any subset $H\subset\mathcal{M}([0,1])$, we have
$$\sup_{\tilde{\mu}\in\hat{\Psi}^{-1}(H)}\widetilde{\mathcal{F}}(\tilde{\mu})
=\sup_{\mu\in H}\mathcal{F}_{\varphi}(\mu).$$
\begin{proof}
First we show the inequality $\displaystyle\sup_{\tilde{\mu}\in\hat{\Psi}^{-1}(H)}\widetilde{\mathcal{F}}(\tilde{\mu})
\le \sup_{\mu\in H}\mathcal{F}_{\varphi}(\mu)$.
Let $\tilde{\mu}\in\hat{\Psi}^{-1}(H)$.
Then it is sufficient to show that
\begin{equation}
\label{com1}
\widetilde{\mathcal{F}}(\tilde{\mu})
\le \sup_{\mu\in H}\mathcal{F}_{\varphi}(\mu).
\end{equation}
If $\tilde{\mu}$ is not $\sigma$-invariant, then $\widetilde{\mathcal{F}}(\tilde{\mu})=-\infty$ and so the
inequality (\ref{com1}) is trivial. Thus, we assume that $\tilde{\mu}$ is $\sigma$-invariant.
Since $\Psi$ is injective away from a countable set, we have $h(\tilde{\mu})=h(\hat{\Psi}(\tilde{\mu}))$.
Thus, we have
\begin{eqnarray*}
\widetilde{\mathcal{F}}(\tilde{\mu})
&=&
h(\tilde{\mu})-\int\varphi\circ\Psi d\tilde{\mu} \\
&=&
h(\hat{\Psi}(\tilde{\mu}))-\int\varphi d\hat{\Psi}(\tilde{\mu}) \\
&=&
\mathcal{F}_{\varphi}(\hat{\Psi}(\tilde{\mu}))\\
&\le&
\sup_{\mu\in H}\mathcal{F}_{\varphi}(\mu).
\end{eqnarray*}
Here the last inequality follows from $\tilde{\mu}\in\hat{\Psi}^{-1}(H)$. Thus, we have (\ref{com1}).

Next, we show the opposite inequality. Let $\mu\in H$.
It is sufficient to show that
\begin{equation}
\label{com2}
\mathcal{F}_{\varphi}(\mu)\le \sup_{\tilde{\mu}\in\hat{\Psi}^{-1}(H)}\widetilde{\mathcal{F}}(\tilde{\mu}).
\end{equation}
Without loss of generality, we may assume that $\mu$ is $f$-invariant.
We set $\tilde{\mu}:=\mu\circ \Phi^{-1}$. Since $\Psi\circ\Phi$ is an identity map,
we have $\hat{\Psi}(\tilde{\mu})=\mu$ and so $\tilde{\mu}\in\hat{\Psi}^{-1}(H)$.
Moreover, since $\Phi$ is injective, we have $h(\mu)=h(\tilde{\mu})$.
Thus, we have
\begin{eqnarray*}
\mathcal{F}_{\varphi}(\mu)
&=&
h(\mu)-\int\varphi d\mu \\
&=&
h(\tilde{\mu})-\int\varphi\circ\Psi d\tilde{\mu} \\
&=&
\widetilde{\mathcal{F}}(\tilde{\mu}),
\end{eqnarray*}
which implies (\ref{com2}).
\end{proof}
\end{lem}

Now, we prove Corollary \ref{p-e}. First we show
$$\liminf_{n\rightarrow\infty}\frac{1}{n}\log L(\{x\in [0,1]:\mathcal{E}_n(x)\in G\})\ge\sup_{\mu\in G}\mathcal{F}_{\varphi}(\mu)$$
for any open set $G\subset\mathcal{M}([0,1])$.
Take any open set $G\subset\mathcal{M}([0,1])$.
By Lemma \ref{conti}, $\hat{\Psi}^{-1}(G)$ is open in $\mathcal{M}(\Sigma_f)$.
Since $(\Sigma_f,\sigma)$ satisfies a level-2 large deviation principle with $L_f$ and the rate function is $\widetilde{\mathcal{F}}$, we have
$$\liminf_{n\rightarrow\infty}\frac{1}{n}\log L_f(\{\omega\in\Sigma_f:\mathcal{E}_n(\omega)\in\hat{\Psi}^{-1}(G) \})
\ge\sup_{\tilde{\mu}\in\hat{\Psi}^{-1}(G)}\widetilde{\mathcal{F}}(\tilde{\mu})
=\sup_{\mu\in G}\mathcal{F}_{\varphi}(\mu),$$
where we use Lemma \ref{com} for the last equality.
By the condition \ref{2.1}, there exists a Borel subset $Y\subset\Sigma_f$ such that
$[0,1]\setminus \Psi(Y)$ is countable and both $\Psi\circ\Phi\colon [0,1]\to [0,1]$ and $\Phi\circ\Psi\colon Y\to \Phi(\Psi(Y))$ are identity.
Then we have
\begin{eqnarray*}
& &
\liminf_{n\rightarrow\infty}\frac{1}{n}\log L(\{x\in [0,1]:\mathcal{E}_n(x)\in G\})\\
&=&
\liminf_{n\rightarrow\infty}\frac{1}{n}\log L(\{x\in \Psi(Y):\mathcal{E}_n(x)\in G\}) \\
&=&
\liminf_{n\rightarrow\infty}\frac{1}{n}\log L_f(\{\omega\in Y:\mathcal{E}_n(\omega)\in \hat{\Psi}^{-1}(G)\}) \\
&=&
\liminf_{n\rightarrow\infty}\frac{1}{n}\log L_f(\{\omega\in\Sigma_f:\mathcal{E}_n(\omega)\in\hat{\Psi}^{-1}(G) \}).
\end{eqnarray*}
Thus we get the lower bound.
The upper bound can be shown similarly.

\section{Applications}
In this section, we apply Theorem \ref{main1} and Corollary \ref{p-e} to systems appeared in \S3.
\subsection{Simple examples}
\begin{prop}
\label{se1ldp}
Let $(I,f)$ be as in Example \ref{se1} and $\varphi\colon I\to \mathbb{R}$ be a constant function $\varphi\equiv \log 3$. Then $(I,f)$ satisfies a level-2 large deviation principle with
the Lebesgue measure $L$ on $[0,1]$, and the rate function is the free energy $\mathcal{F}_{\varphi}\colon\mathcal{M}([0,1])\to [-\infty,0]$ given by (\ref{free}).
\begin{proof}
Let $(\Sigma_f,\sigma)$ be the symbolic space of $(I,f)$.
We set $J_0=[0,\frac{1}{6})$, $J_1=[\frac{1}{6},\frac{1}{2})$, $J_2=[\frac{1}{2},\frac{2}{3})$, $J_3=[\frac{2}{3},\frac{5}{6})$ and $J_4=[\frac{5}{6},1]$.
Define a measurable map $\Phi\colon I\to\Sigma_f$ by
$$(\Phi(x))_i=j\text{ if and only if }f^i(x)\in J_j.$$
Then it is not difficult to see that $\Phi$ satisfies \ref{2.1}.

Thus, it is sufficient to show that $(\Sigma_f,\sigma)$ satisfies the condition \ref{2.2} of Corollary \ref{p-e}
for $L_f=L\circ\Phi^{-1}$ and the constant function $\varphi\circ\Psi\equiv\log 3$.
As we have shown in Example \ref{se1}, $(\Sigma_f,\sigma)$ satisfies
the one-way specification property, and so  \ref{1.1} holds.
Since the slope of $f$ is $3$, and the definition of $\Phi$, we have
$$\left(\frac{1}{2}\right)3^{-n}\le L\circ\Phi^{-1}([w])\le 3^{-n}.$$
This implies that \ref{1.2} and \ref{1.3} hold for $L_f$ and $\varphi\circ \Psi$.
Thus we have \ref{2.2}
and so the proposition is proved.
\end{proof}
\end{prop}

\begin{prop}
\label{se2ldp}
Let $(S^1,f)$ be as in Example \ref{se2} and $\varphi\colon S^1\to \mathbb{R}$ be a continuous function
defined by $\varphi:=\max\{\log f',0\}$. 
Then $(S^1,f)$ satisfies a level-2 large deviation principle with
the Lebesgue measure $L$ on $S^1$ and the rate function is the free energy $\mathcal{F}_{\varphi}\colon\mathcal{M}(S^1)\to [-\infty,0]$ given by (\ref{free}).
\begin{proof}
In what follows we will show that $(S^1,f)$, $L$, and $\varphi$ satisfy the conditions \ref{A.I}-\ref{A.IIII} of Theorem \ref{main1}.
$(S^1,f)$ satisfies the one-way specification property with $\{0\}$ and $\{\frac{1}{2}\}$ (see Example \ref{se2}) and so \ref{A.I} holds.
Since $\Omega(f)=\{0,\frac{1}{2}\}$, we have $h(\mu)=0$ for any $\mu\in\mathcal{M}_f(S^1)$, and so \ref{A.III} holds.
Since $f$ is a one-dimensional map, by the standard mean value theorem, we can see that $L$ and $\varphi$ satisfies \ref{A.II} and \ref{A.IIII}.
Thus the proposition is proved.
\end{proof}
\end{prop}

\subsection{($-\beta$)-transformations; Proof of Theorem \ref{minus beta}}
The aim of this subsection is to give a proof of Theorem \ref{minus beta}. Let $\beta>1$ be a Yrrap number.
First we define a map $\Phi\colon [0,1]\to \Sigma_{-\beta}$ as follows:
If $i'(1)=(\cdots 0)^{\infty}$, then let $J_0=[0,\frac{1}{\beta}]$, $J_i=(\frac{i}{\beta},\frac{i+1}{\beta}]$ for $1\le i\le b-1$ and
$J_b=(\frac{b}{\beta},1]$.
If $i'(1)$ is not of the form $(\cdots 0)^{\infty}$, then let
$J_i=[\frac{i}{\beta},\frac{i+1}{\beta})$ for $0\le i\le b-1$ and $J_b=[\frac{b}{\beta},1]$.
Then, we define $\Phi\colon [0,1]\to \Sigma_{-\beta}$ by
$(\Phi(x))_i=j$ if and only if $T_{-\beta}(x)\in J_j$.
Then, it is not difficult to see that $\Phi$ satisfies \ref{2.1}.

To prove Theorem \ref{minus beta}, it is sufficient to show that $(\Sigma_{-\beta},\sigma)$ satisfies the
condition \ref{2.2} for $L_{-\beta}=L\circ\Phi^{-1}$ and the constant function $\varphi\circ\Psi\equiv\log\beta$. 
By Proposition \ref{beta-OWSP}, $(\Sigma_{-\beta},\sigma)$ satisfies \ref{1.1}.
Moreover, for $w\in\mathcal{L}_n(\Sigma_{-\beta})$, it is easy to see that $\Phi^{-1}([w])$ is a sub-interval of
$(0,1]$ and $T_{-\beta}^n\colon \Phi^{-1}([w])\to (0,1]$ is a
continuous linear map whose slope is equal to $\beta^n$.
Thus we have
$$L_{-\beta}([w])=L(\Phi^{-1}([w]))\le\beta^{-n}.$$
This implies \ref{1.3}. In what follows we will show that the condition \ref{1.2} holds for
$L_{-\beta}$ and the constant function $\log \beta$.
\begin{lem}
\label{l-cylinder}
Let $w\in\mathcal{L}_n(\Sigma_{-\beta})$ and
suppose that there exist $0\le c<d\le b$ such that
$wc,wd\in\mathcal{L}(\Sigma_{-\beta})$. Then
\begin{equation}
\label{l-cylinder0}
L_{-\beta}([w])\ge \left(1-\frac{b}{\beta}\right)\beta^{-n}.
\end{equation}
\begin{proof}
It is easy to see that if $(a_1,a_2)\subset J_i$ and there exists a point
$b_1\in J_j$ such that $T_{-\beta}(b_1)=a_1$, then there exists a point
$b_2\in J_j$ such that $T_{-\beta}(b_2)=a_2$ and $\beta(b_1-b_2)=a_2-a_1$.
Thus, for any $x\in \Phi^{-1}([w])$,
$$T_{-\beta}^{-1}(\cdots(T_{-\beta}^{-1}(T_{-\beta}^{-1}(w_n\beta^{-1},
T_{-\beta}^{n-1}(x))\cap J_{w_{n-1}})\cap J_{w_{n-2}})\cap\cdots)\cap J_{w_1})$$
is an interval of length $\beta^{-n+1}|T_{-\beta}^{n-1}(x)-w_n\beta^{-1}|$
and is contained in $\Phi^{-1}([w])$. So we have
\begin{equation}
\label{l-cylinder1}
L_{-\beta}([w])\ge \beta^{-n+1}|T_{-\beta}^{n-1}-w_n\beta^{-1}|
\end{equation}
for any $x\in \Phi^{-1}([w])$.
On the other hand, $T_{-\beta}(w_n\beta^{-1})=1$ and so $wb\in\Sigma_{-\beta}$.
By the assumption, there exists a point $x\in \Phi^{-1}([w])$ such that
$T_{-\beta}^n(x)\not\in J_b$.
Then we have
$$|T_{-\beta}^n(x)-T_{-\beta}(w_n\beta^{-1})|\ge |J_b|=\left(1-\frac{b}{\beta}\right),$$
which implies
\begin{equation}
\label{l-cylinder2}
|T_{-\beta}^{n-1}(x)-w_n\beta^{-1}|\ge\beta^{-1}\left(1-\frac{b}{\beta}\right).
\end{equation}
Combining the equations (\ref{l-cylinder1}) and (\ref{l-cylinder2}),
we have (\ref{l-cylinder0}).
\end{proof}
\end{lem}

For $w\in\mathcal{L}(\Sigma_{-\beta})$,
let $g_{\beta}(w)$ be a minimum integer $i\ge 0$ so that there exist $v\in\mathcal{L}_i(\Sigma_{-\beta})$ and
integers $0\le c<d\le b$ such that $wvc,wvd\in\mathcal{L}(\Sigma_{-\beta})$,
and set
$$g_{\beta}(n):=\sup_{w\in\mathcal{L}_n(\Sigma_{-\beta})}g_{\beta}(w).$$
for $n\ge 1$.

\begin{lem}
\label{l-energy}
Suppose that $\displaystyle\lim_{n\rightarrow\infty}\frac{g_{\beta}(n)}{n}=0$.
Then the condition \ref{1.2} holds for $L_{-\beta}$ and the constant function $\log \beta$.
\begin{proof}
Let $\epsilon>0$. By the assumption, there exists an integer $N$
such that for any $n\ge N$, $g_{\beta}(n)\le n\epsilon$ holds.
Fix any $n\ge N$ and $w\in\mathcal{L}_n(\Sigma_{-\beta})$.
Then we have $g_{\beta}(w)\le g_{\beta}(n)\le n\epsilon$.
Thus we can find $v\in\mathcal{L}(\Sigma_{-\beta})$ and
$0\le c<d\le b$ such that
$$|v|\le n\epsilon\text{ and }wvc,wvd\in\mathcal{L}(\Sigma_{-\beta}) $$
holds. So it follows from Lemma \ref{l-cylinder} that
\begin{eqnarray*}
L_{-\beta}([w])
&\ge &
L_{-\beta}([wv]) \\
&\ge &
\left(1-\frac{b}{\beta}\right)\beta^{-|wv|} \\
&\ge &
\left(1-\frac{b}{\beta}\right)\beta^{-n(1+\epsilon)}.
\end{eqnarray*}
Since $w\in\mathcal{L}_n(\Sigma_{-\beta})$, $n\ge N$ and $\epsilon>0$
are arbitrary, the above inequality implies that
\ref{1.2} holds for $L_{-\beta}$ and the constant function $\log \beta$.
\end{proof}
\end{lem}

Thus, it follows from Lemmas \ref{beta-lem1} and \ref{l-energy} that the condition \ref{1.2} holds for $L_{-\beta}$ and the constant function $\log\beta$.
Therefore, $(\Sigma_{-\beta},\sigma)$ satisfies the condition \ref{2.2} for $L_{-\beta}$ and the constant function $\varphi\circ\Psi\equiv\log\beta$,
and so Theorem \ref{minus beta} is proved.

\begin{ex}
\label{minimal pisot}
Finally, we give an example of the $(-\beta)$-transformation whose rate functions with respect to
Lebesgue measure and maximal entropy measure do not coincide.
Let $\beta$ be the minimal Pisot number, that is, the real root of $X^3-X-1=0$.
Then it is easy to see that $i'(1)=100(1)^{\infty}$ and so, $\beta$ is a Yrrap number and $(\Sigma_{-\beta},\sigma)$
is generated by a graph $\Gamma_{-\beta}$, which is shown in Figure \ref{fig:5}.
Thus, we have
\begin{equation}
\label{6pisot}
\mathcal{M}_{\sigma}(\Sigma_{-\beta})=\{\delta_{0^{\infty}}\}\cup \{\delta_{1^{\infty}}\}\cup\mathcal{M}_{\sigma}(X_3).
\end{equation}
Here $X_3$ is a subshift generated by a graph $\Gamma_{2,\infty}$ shown in Figure \ref{fig:6}.
In particular, we can see that $\mathcal{M}_{\sigma}^e(\Sigma_{-\beta})$ is not dense in $\mathcal{M}_{\sigma}(\Sigma_{-\beta})$.
It follows from \cite[Corollary 2.4]{LS} that $T_{-\beta}$ (resp. $(\Sigma_{-\beta},\sigma)$) admits a unique measure of maximal entropy
$m$ (resp. $m_{-\beta}$).
Then by (\ref{6pisot}), $m_{-\beta}$ is supported on $X_3$.
Since $(X_3,\sigma)$ satisfies the (W)-specification,
by \cite[Theorem A]{CTY} with an easy additional calculation, we can show that
$(\Sigma_{-\beta},\sigma)$ satisfies a level-2 large deviation principle with $m_{-\beta}$ and the rate function
$q^{m_{-\beta}}\colon\mathcal{M}(\Sigma_{-\beta})\to [-\infty,0]$ is given by
$$q^{m_{-\beta}}(\mu)=
\left\{
\begin{array}{ll}
h(\mu)-\log\beta & (\mu\in\mathcal{M}_{\sigma}(X_3)); \\
-\infty & (\text{otherwise}).
\end{array}
\right.
$$
So in a similar way to the proof of Corollary \ref{p-e}, we can show that
$([0,1],T_{-\beta})$ satisfies a level-2 large deviation principle with $m$ and
the rate function $q^m\colon\mathcal{M}([0,1])\to[-\infty,0]$ is given by
$$q^m(\mu)=
\left\{
\begin{array}{ll}
h(\mu)-\log\beta & (\mu\in \hat{\Psi}^{-1}(\mathcal{M}_{\sigma}(X_3))); \\
-\infty & (\text{otherwise}).
\end{array}
\right.
$$
On the other hand, by Theorem \ref{minus beta}, $([0,1],T_{-\beta})$ satisfies a level-2 large deviation principle with the Lebesgue measure
and the rate function $q^L\colon\mathcal{M}([0,1])\to [-\infty,0]$ coincides with the free energy $\mathcal{F}_{\varphi}$ given by (\ref{free}) with respect to a constant function
$\varphi\equiv \log \beta$.
Thus, we have $q^m(\mu)\not =q^L(\mu)$ for $\mu\in\mathcal{M}_{T_{-\beta}}([0,1])\setminus\hat{\Psi}^{-1}(\mathcal{M}_{\sigma}(X_3))$.
Moreover, it is not difficult to see that
$$\mathcal{M}_{T_{-\beta}}([0,1])\setminus\hat{\Psi}^{-1}(\mathcal{M}_{\sigma}(X_3))=\left\{a\delta_{\frac{1}{\beta+1}}+(1-a)\mu:0<a\le 1,\mu\in\mathcal{M}_{T_{-\beta}}((0,1])\right\}.$$
In particular, this set is not empty. Thus, we conclude that $q^m\not=q^L$.
\end{ex}

\vspace*{3mm}

\noindent
\textbf{Acknowledgement.}~ 
This work was supported by JSPS KAKENHI Grant Number JP16K17611.

\end{document}